\documentclass[12pt,a4paper]{amsart}
\usepackage{amsfonts}
\usepackage{amsthm}
\usepackage{dsfont}

\usepackage{amsmath}
\usepackage{amscd}
\usepackage[latin2]{inputenc}
\usepackage{t1enc}
\usepackage[mathscr]{eucal}
\usepackage{indentfirst}
\usepackage{graphicx}
\usepackage{graphics}
\usepackage{pict2e}
\usepackage{epic}
\numberwithin{equation}{section}
\usepackage[margin=2.9cm]{geometry}
\usepackage{epstopdf} 
\DeclareMathOperator\supp{supp}

\newcommand{\of}[1]{\left(#1\right)}
\usepackage{xcolor}
 
\newcommand{\floor}[1]{\lfloor #1 \rfloor}

\theoremstyle{plain}
\newtheorem{thm}{Theorem}[section]
\newtheorem{lemma}[thm]{Lemma}

\theoremstyle{definition}

\newtheorem{definition}[thm]{Definition}

\theoremstyle{remark}

\newtheorem*{stat*}{Statement}
\newtheorem*{prop*}{Proposition}
\numberwithin{equation}{section}
\newcommand{\E}{\mathbb{E}}

\newcommand{\N}{\mathbb{N}}

\newcommand{\R}{\mathbb{R}}

\usepackage{hyperref}
\DeclareMathOperator{\var}{Var}
 \hypersetup{
    colorlinks=true,
    linkcolor=blue,
    urlcolor=blue,
    linktoc=all
}

\begin{document}

\title{A note on pair dependent linear statistics with slowly growing variance.}

\author[A. Aguirre and A. Soshnikov]{Ander Aguirre and Alexander Soshnikov}

\address{University of California at Davis \\ Department of Mathematics \\ 1 Shields Avenue \\Davis CA 95616 \\ United States of America} 

\email{aaguirre@ucdavis.edu}

\address{University of California at Davis \\ Department of Mathematics \\ 1 Shields Avenue \\  Davis CA 95616 \\ United States of America} 

\email{soshniko@math.ucdavis.edu}

\begin{abstract} 
We prove Gaussian fluctuation for
pair counting statistics of the form \\
$ \sum_{1\leq i\neq j\leq N} f(\theta_i-\theta_j)$
for the Circular Unitary Ensemble (CUE) of 
random matrices in the case of a slowly growing variance in the limit of large $N.$
\end{abstract}

\maketitle

\section{Introduction} 
\vspace{5mm}

The study of random matrix theory goes back to the principal component data analysis of J. Wishart in 1920s-1930s and revolutionary ideas of E.Wigner in quantum physics in 1950s that linked
statistical properties of the energy levels of heavy-nuclei atoms with spectral properties of Hermitian random matrices with independent components.

In 1960s, F. Dyson introduced three archetypal types of matrix ensembles: Circular Orthogonal Ensemble (COE), Circular Unitary Ensemble (CUE),
and Circular Symplectic Ensemble (CSE), see e.g. \cite{Dyson1}-\cite{Dyson4}. 
The probability density of the eigenvalues $ \{e^{i\*\theta_j}\}_{j=1}^N, \ \ 0\leq \theta_1, \ldots, \theta_N <2\*\pi,$ is given by
 \begin{align}
            \label{betaensemble}
                p_N(\overline{\theta})=\frac{1}{Z_{N}(\beta)}\prod_{1\leq j< k\leq N}\left|e^{i\theta_j}-e^{i\theta_k}\right|^{\beta}, 
            \end{align}
where $\beta=1,2,$and $4$ correspond to COE, CUE, and CSE, correspondingly.
For arbitrary $\beta>0,$ a (sparse) random matrix model with eigenvalues distribution following (\ref{betaensemble}) was introduced in \cite{KN}. The ensemble (\ref{betaensemble}) for arbitrary $\beta>0$ is known as the Circular Beta Ensemble (C$\beta$U).

 The Circular Unitary Ensemble ($\beta=2$) corresponds to the joint distribution of the eigenvalues of an $N\times N$ random unitary matrix U distributed according to the Haar measure. In
particular, the partition function is given by
\begin{align}
\label{partfun}
Z_{N}(2)=(2\pi)^N \times N!
\end{align}

In \cite{paper}, pair counting statistics of the form
            \begin{align}
\label{pairs}
                S_N(f)=\sum_{1\leq i\neq j\leq N} f(L_N\*(\theta_i-\theta_j)),
            \end{align}
were studied for C$\beta$E (\ref{betaensemble}) for $1\leq L_N\leq N$ under certain smoothness assumptions on $f.$ The research in \cite{paper} was motivated by a classical result of Montgomery on pair correlation of zeros of the Riemann zeta function \cite{montgomery1}-\cite{montgomery2}.
Assuming the Riemann Hypothesis, Montgomery studied the distribution of the 
``non-trivial'' zeros on the critical line $1/2 +i\*\R.$ Rescaling zeros $\{ 1/2 \pm \gamma_n\},$
\[ 
\tilde{\gamma}_n= \frac{\gamma_n}{2\*\pi} \*\log(\gamma_n),
\]
Montgomery considered the statistic $\sum_{0<\tilde{\gamma}_j\neq \tilde{\gamma}_k <T} f(\tilde{\gamma}_j- \tilde{\gamma}_k )$ for large $T$ and sufficiently fast decaying $f$ with $\supp{\mathcal{F}(f)}  \subset [-\pi, \pi],$ where $\mathcal{F}(f)$ denotes the Fourier transform of $f.$ The results of \cite{montgomery1}-\cite{montgomery2} imply that the two-point correlations of the (rescaled) critical zeros coincide in the limit with the local two point correlations of the eigenvalues of a CUE random matrix.

The results of \cite{paper} deal with the limiting behavior of (\ref{pairs}) in three different regimes, namely macroscopic ($L_N=1$), mesoscopic ($1\ll L_N \ll N$) and microscopic ($L_N=N$).  
In the unscaled $L_N=1$ case it was shown that
    \begin{align}
\label{pairs1}
S_N(f)=\sum_{1\leq i\neq j\leq N} f(\theta_i-\theta_j),
            \end{align}
has a non-Gaussian fluctuation in the limit $N\to \infty$ provided $f$ is a sufficiently smooth function on the unit circle.

Namely, let $f$ be a real even integrable function on the unit circle.
Denote the Fourier coefficients of $f$ as
\begin{align}
\label{FourierS}
\hat{f}(k)=\frac{1}{2\*\pi}\*\int_0^{2\*\pi}f(x)\* e^{-i\*k\*x}\* dx.
\end{align}
Let us assume that $f'\in L^2(\mathbb{T})$ for $\beta=2$,
$\sum_{k\in \mathbb{Z}}|\hat{f}(k)|\*|k|<\infty$ for $\beta<2, \ \sum_{k\in \mathbb{Z}}|\hat{f}(k)|\*|k|\*\log(|k|+1)<\infty$ for $\beta=4,$ 
and 
$\sum_{k\in \mathbb{Z}}|\hat{f}(k)||k|^2<\infty$ for $\beta \in (2,4)\cup (4, \infty).$
Then we have the following convergence in distribution as $N\rightarrow \infty$:

            $$S_N(f)-\E S_N(f)\xrightarrow{\hspace{2mm}\mathcal{D}\hspace{2mm}  } \frac{4}{\beta}\sum_{k=1}^{\infty}\hat{f}(k)k(\varphi_k-1),$$
        
        where $\varphi_m$ are i.i.d. exponential random variables with  $\E(\varphi_m)=1$. 
For $\beta=2$ the result was proven under the optimal condition $\sum_{k\in \mathbb{Z}}|\hat{f}(k)|^2\*|k|^2<\infty.$ 

The goal of this paper is to study the fluctuation of the pair counting statistic (\ref{pairs1}) when $\var(S_N(f))$ slowly grows with $N$ to infinity.

\begin{definition} A positive sequence ${V_N}$ is said to be slowly varying in sense of Karamata (\cite{BGT}) if

\begin{align}
    \label{karamata}
\lim_{n\rightarrow \infty}\frac{V_{\floor{\lambda N}}}{V_N}= 1, \quad \quad \forall\lambda>0,
\end{align}
where $\floor{m}$ denotes the integer part of $m$ .

\end{definition}
The following notation will be used throughout the paper:
\begin{align}
\label{vvvv}
V_N=\sum_{k=-N}^{k=N}|\hat{f(k)}|^2|k|^2.  
\end{align}
\begin{thm} Let $f\in L^2(\mathbb{T})$ be a real even function such that $V_N=\sum_{k=-N}^{k=N}|\hat{f(k)}|^2|k|^2, \\N=1,2,\ldots, $ is a slowly varying sequence that diverges to infinity as $N\to \infty$. Then we have the following convergence in distribution
 \[
\frac{ S_N(f)-\E  S_N(f)}{\sqrt{2\*\sum_{-N}^{N} |\hat{f}(k)|^2\*|k|^2}}\xrightarrow{\hspace{2mm}\mathcal{D}\hspace{2mm}  }
\mathcal{N}(0,1)
\]
\end{thm}

Linear statistics of the eigenvalues of random matrices $\sum_{j=1}^N f(\lambda_j) $ have been studied extensively in the literature. Johansson (\cite{johansson1}) proved for (\ref{betaensemble}) for arbitrary $\beta>0$ and sufficiently smooth real-valued $f$ that
\[
\frac{\sum_{j=1}^N f(\theta_j) - N\*\hat{f}(0)}{\sqrt{\frac{2}{\beta}\*\sum_{-\infty}^{\infty} |\hat{f}(k)|^2\*|k|}}
\]converges in distribution to a standard Gaussian random variable. In particular, for $\beta=2$ he proved the result under the optimal conditions on $f$, namely 
\[\sum_{-\infty}^{\infty} |\hat{f}(k)|^2\*|k|<\infty.\]
If the variance of the linear statistic goes to infinity with $N,$ Diaconis and Evans \cite{DE} 
proved the CLT in the case 
$\beta=2$ provided the sequence $\{\sum_{-N}^{N} |\hat{f}(m)|^2\*|m| \}_{n \in \N}$ is slowly varying.

For the results on the linear eigenvalue statistics in the mesoscopic regime\\
$\sum_{j=1}^N f(L_N\*\theta_j), \ 1\ll L_N \ll N,$ we refer the reader to \cite{sasha}, \cite{BL}, \cite{lambert}, and references therein. For additional results on the spectral properties of C$\beta$U we refer the reader to \cite{DS}, \cite{johansson3}, \cite{BF}, \cite{Sasha}, \cite{HKOC}, \cite{meckes}, \cite{PZ}, \cite{webb}, \cite{WF1}, and references therein.

The proof of the main result of the paper (Theorem 1.2) is given in the next section. Throughout the paper, he notation $a_N=O(b_N)$ means 
that the ratio $ a_N/b_N$ is bounded from above in absolute value. The notation $a_N=o(b_N)$ means that $a_n/b_N\to 0$ as $N\to \infty.$
Occasionally, for non-negative quantities, in this case we will also use the notation $a_N \ll b_N.$

Research has been partially supported  by the Simons Foundation Collaboration Grant for Mathematicians \#312391.

\vspace{5mm}

\section{Proof of Theorem 1.2}

\vspace{5mm}
The section is devoted to the proof of Theorem 1.2.
We start by recalling the formula for the variance of $S_N(f)$ from Proposition 4.1 of \cite{paper}: 
\begin{align}
\label{v1}
            \var(S_N(f)) &=4\*\sum_{1\leq s\leq N-1}s^2|\hat{f}(s)|^2 + 4\*
(N^2-N)\sum_{N\leq s} |\hat{f}(s)|^2 \\
            &-4\sum_{\substack{1\leq s,t \\ 1\leq |s-t|\leq N-1\\ N\leq \max(s,t)}}(N-|s-t|)\hat{f}(s)\hat{f}(t)\hspace{2mm}
            -4\sum_{\substack{1\leq s,t\leq N-1\\N+1\leq s+t}}((s+t)-N) \hat{f}(s)\hat{f}(t). \nonumber
        \end{align}

Our first goal is to show that the last two (off-diagonal) terms in the variance expression (\ref{v1}) are much smaller than $V_N =\sum_{s=-N}^{N} s^2|\hat{f}(s)|^2$ for large $N$ provided (\ref{karamata})
is satisfied.

  \begin{lemma} Let  $V_N$ from (\ref{vvvv}) be a slowly varying sequence diverging to infinity as $N\to \infty$. Then, as $N\to\infty$, we have
            \begin{enumerate}
                \item[(i)]
                    \begin{align*}
                        \sum_{\substack{1\leq s,t\leq N\\ s+t\geq N+1}}s|\hat{f}(s)|\cdot|\hat{f}(t)| =o\of{V_N},
                    \end{align*}
                \item[(ii)]
                    \begin{align*}
                        (N+1)\sum_{\substack{s-t\leq N\\s\geq N+1\\1\leq t\leq N}}|\hat{f}(s)|\cdot|\hat{f}(t)|=o\of{V_N} ,
                    \end{align*}    
                \item[(iii)]
                    \begin{align*}
                        N\sum_{\substack{|s-t|\leq N-1\\s,t\geq N}}|\hat{f}(s)|\cdot|\hat{f}(t)|\ =o\of{V_N}.
                    \end{align*}
            \end{enumerate}
            
        \end{lemma}
The proof of the lemma is somewhat similar to the proof of Lemma 4.4. in \cite{paper}.  For the convenience of the reader, we give the full details of the proof below.

\begin{pfo}{\textit{Lemma 2.1}}

 Proof of (i).  

Let $x_s=s|\hat{f}(s)|$ for $1\leq s \leq N$ and $X_N=\{x_s\}_{s=1}^N$. 
Define a vector $Y_N:=X_N\*\mathds{1}_{(s> N/2)},$  so that the first $\floor{N/2}$ 
coordinates of $Y_N$ are zero and the rest coincide with the corresponding coordinates 
of $X_N.$  We note that 
\begin{align}
\label{tom}
2\*||X_N||_2^2=V_N \ \text{and} \ 
||Y_N||_2^2=o(V_N),
\end{align}
where $||X||_2$ denotes the Euclidean norm of a vector $X\in \R^N.$
The last bound follows from the condition (\ref{karamata}) on the slow growth of $V_N.$

We now write the off diagonal variance term in (i) as a bilinear form: 
                    \begin{align}
                    \label{4.15}
                        \sum_{\substack{1\leq s,t\leq N\\ s+t\geq N+1}}s|\hat{f}(s)|\cdot|\hat{f}(t)|
                        &=\sum_{t=1}^{N} x_t \cdot \left(\frac{1}{t}\sum_{s=N-t+1}^{N} x_s\right)\nonumber\\
                        &=\sum_{t=1}^N x_t \cdot \left(\frac{1}{t}\sum_{s=1}^t (U_N\*X_N)_s\right)
                        =\langle X_N,A_N X_N\rangle,
                    \end{align}
with $A_N=B_N\*U_N$, where $U_N$ is a unitary permutation 
matrix given by $(U_N)_{s,t}= \mathds{1}_{(t=N-s+1)}$ and $B_N$ is a lower triangular matrix given by $(B_N)_{s,t}= (1/s)\mathds{1}_{(t\leq s)}$. 
The matrix $A_N$ is given by:
            \[
                A_N = 
                        \renewcommand\arraystretch{1.25}
                        \begin{pmatrix}
                            0           & 0           & 0           & 0  & \dots  & 1 \\
                              0       &  0  &\ddots &\ddots   &\frac{1}{2} & \frac{1}{2}  \\
                             0 & \ddots & \ddots  &\frac{1}{3} & \frac{1}{3} & \frac{1}{3} & \\
                            \vdots & \vdots & \vdots  & \vdots      & \vdots      & \vdots \\
                            \frac{1}{N} & \frac{1}{N} &\frac{1}{N}  & \frac{1}{N} &\dots  & \frac{1}{N}
                        \end{pmatrix}
           \]
Taking into account that the upper-right $\floor{N/2}\times \floor{N/2}$ block of $A_N$ is zero, we can bound
$\langle X_N , A_N X_N \rangle\leq \langle X_N, A_N Y_N\rangle + \langle Y_N , A_N X_N \rangle. $   It was shown in \cite{paper} that $||A_N||_{op}\leq 3,$ where $||A||_{op}$ denotes the operator norm.
This implies that the expression in (\ref{4.15})is bounded from above by $3\*||X_N||_2\*||Y_N||_2=o(V_N)$  by (\ref{tom}). This completes the proof of Lemma 2.1(i).\\

To prove part (ii), let $B_N$ be defined as in the proof of part $(i)$. Similarly, 
let $x_s=s|\hat{f}(s)|, \ 1\leq s\leq 2\*N,$ and $X_N=\{x_s\}_{s=1}^{2N}$. 
Now, $X_N$ is a $2N$-dimensional vector such that $||X_N||_2^2$  grows slowly in $N$. 
Define $Y_N:=X_N\*\mathds{1}_{(s> N)},$  so that the first $N$ 
coordinates of $Y_N$ are zero and the rest coincide with the corresponding coordinates 
of $X_N.$ 
Observe that
                \begin{align*}
                    N\sum_{\substack{s-t\leq N\\s\geq N+1\\1\leq t\leq N}}|\hat{f}(s)|\cdot|\hat{f}(t)|
                    &\leq \sum_{t=1}^N x_t\left(\frac{1}{t}\sum_{s=N+1}^{N+t}x_s\right)\\
                    &= \langle  C_N X_N, D_N X_N\rangle,
                \end{align*}
            where 
               
                 $$ C_N=\renewcommand\arraystretch{1.25}
                        \begin{pmatrix}
                           I_N & 0\\
                           0   & 0
                        \end{pmatrix}
                        \hspace{5mm}\text{and}\hspace{5mm}
                    D_{N} = 
                        \renewcommand\arraystretch{1.25}
                        \begin{pmatrix}
                           0 & B_N\\
                           0 & 0
                        \end{pmatrix}$$
We note that $||D_N||_{op}\leq 3$ and  $||C_N||_{op}=1.$  Once again:
$$\langle  C_N X_N, D_N X_N\rangle=
\langle  C_N X_N, D_N Y_N\rangle \leq 3\*||X_N||_2\*||Y_N||_2=o(V_N).$$
The proof of (ii) is completed.\\

Proof of (iii).
We start by estimating the l.h.s. of (iii) from above by
                \begin{align}
\label{myt}
                 2\*N\*\sum_{\substack{t-N+1
                 \leq s\leq N+t-1\\t\geq s\geq N}}|\hat{f}(s)|\cdot|\hat{f}(t)|.
               \end{align}
                As before, let $x_s=s|\hat{f}(s)|$ for $s\geq 1.$ 
Define $X=\{x_s\}_{s=1}^{\infty}$ and \\
$X^{(j)}=X\*\mathds{1}_{(j\*N\leq s< (j+1)\*N)}, \ j=0,1,2,\ldots.$
We can bound (\ref{myt}) from above by 
                    \begin{align}
                    \label{jj}
                    2\* \sum_{t=N}^\infty x_t\left(\frac{1}{t}\sum_{s=t-N+1}^{t}x_s\right)
                        = 2\*\sum_{j=1}^{\infty} \sum_{t=j\*N}^{(j+1)\*N-1} x_t\left(\frac{1}{t}\sum_{s=t-N+1}^{t}x_s\right),
\end{align}
One can write the second sum at the r.h.s. of (\ref{jj}) as 
\begin{align}
\label{jjj}
\sum_{t=j\*N}^{(j+1)\*N-1} x_t\left(\frac{1}{t}\sum_{s=t-N+1}^{t}x_s\right)=
\langle X^{(j)}, R_{N,j} (X^{(j-1)}+X^{(j)}) \rangle, 
\end{align}
where $R_{N,j}$ is a bounded linear operator such that 
\begin{align}
\label{jjjj}
(R_{N,j})_{t,s}=\frac{1}{t}\*\mathds{1}_{(t-N+1\leq s \leq t)}\*\mathds{1}_{(jN\leq t < (j+1)N)}.
\end{align}
The operator norm of $R_{N,j}$ is bounded from above by its Hilbert-Schmidt norm 
                \begin{align}
                \label{j5}
                    ||R_{N,j}||_{op}\leq ||R_{N,j}||_2 = \sqrt{N\*\sum_{t=jN}^{(j+1)N-1} \frac{1}{t^2}}\leq \sqrt{\frac{N^2}{j^2\*N^2}}=\frac{1}{j}.
                \end{align}
            Now, by the Cauchy-Schwarz inequality, the r.h.s. of (\ref{jjj}) can be bounded from above by
                \begin{align}
                \label{j6}
               &\langle X^{(j)}, R_{N,j} (X^{(j-1)}+X^{(j)}) \rangle  \leq ||X^{(j)}||_2\*||R_{N,j}||_{op} \* (||X^{(j-1)}||_2 + ||X^{(j)}||_2)\\
               \label{j7}
               &\leq \frac{1}{j}\*||X^{(j)}||_2^2 +\frac{1}{j}\*||X^{(j)}||_2\* ||X^{(j-1)}||_2.
               \end{align}
Summing up the r.h.s. of the last inequality with respect to $j\geq 1$ gives $o(V_N).$
Indeed, $2\*||X^{(0)}||_2^2=V_N$ and summation by parts gives
\begin{align}
\label{sumparts}
\sum_{j=1}^{\infty}\frac{1}{j}\*||X^{(j)}||_2^2 \leq \sum_{j=1}^{\infty} \frac{1}{j^2}\*(V_{j\*N}-V_N).
\end{align}
The condition (\ref{karamata}) on the slow growth of $V_N$ (\ref{karamata}) implies that the r.h.s. of
(\ref{sumparts}) is $o(V_N)$.
To sum the second term in (\ref{j7}), we write
\begin{align}
\label{0330}
\sum_{j=1}^{\infty} \frac{1}{j}\*||X^{(j)}||_2\* ||X^{(j-1)}||_2= ||X^{(1)}||_2\* ||X^{(0)}||_2 +
\sum_{j=2}^{\infty} \frac{1}{j}\*||X^{(j)}||_2\* ||X^{(j-1)}||_2.
\end{align}
The first term on the r.h.s. of (\ref{0330}) is $\sqrt{(V_{2N}-V_{N})\*V_N}=o(V_N).$
To deal with the second term, we use the Cauchy-Schwarz inequality and proceed as in (\ref{sumparts}).
This completes the proof of Lemma 2.1.
            \end{pfo}
            
\begin{pfo}{\textit{Theorem 1.2}}\\
We now proceed to finish the proof of Theorem 1.2.
We recall that
\begin{align}
\label{formula}
            S_N(f)=\sum_{1\leq i\neq j\leq N}f(\theta_i-\theta_j)=2\sum_{k=1}^{\infty}\hat{f}(k)\left|\sum_{j=1}^N\exp\of{ik\theta_j}\right|^2+ 
\hat{f}(0)\*N^2-N\*f(0).
        \end{align}
Denote by $t_{N,k}$ the trace of the $k$-th power of a CUE matrix, i.e.
\begin{align}
\label{kkk}
t_{N,k}:=\sum_{j=1}^N e^{i\*k\*\theta_j}, \ \ k=0,\pm 1, \pm 2, \ldots.
\end{align}
Then
\begin{align}
&\frac{S_N(f)-\E S_N(f)}{\sqrt{V_N}}=\frac{2}{\sqrt{V_N}}\*\sum_{k=1}^{\infty}\hat{f}(k)\*(|t_{N,k}|^2-\E |t_{N,k}|^2) =\nonumber \\
\label{aaa}  
&\frac{2}{\sqrt{V_N}}\*\sum_{k=1}^{\floor{N/M_N}}\hat{f}(k)\*(|t_{N,k}|^2-\E |t_{N,k}|^2) +\frac{2}{\sqrt{V_N}}\sum_{\floor{N/M_N}+1}^{\infty}\hat{f}(k)\*(|t_{N,k}|^2-\E |t_{N,k}|^2).
\end{align}
Here $\{M_N\}_{N=1}^{\infty}$ is positive integer-valued sequence sufficiently slowly growing to infinity as $N\to \infty$ in such a way that 
\begin{align}
    \label{Kara}
    \lim_{N\to \infty} \frac{V_{\floor{N\*M_N}}}{V_N}=\lim_{N\to \infty} \frac{V_{N}}{V_{\floor{N/M_N}}}=1.
\end{align}
The existence of such a sequence follows from (\ref{karamata}).

It follows from Lemma 2.1 that the second moment of the second term on the r.h.s. of (\ref{aaa}) is
going to zero as $N \to \infty.$ We formulate this result in the next lemma.

\begin{lemma} 
\begin{align}
\label{lemma22}
\var \left(\sum_{\floor{N/M_N}+1}^{\infty}\hat{f}(k)\*|t_{N,k}|^2\right)=o(V_N).
\end{align}

\end{lemma}
\begin{pfo}{\textit{Lemma 2.2}}

It follows from (\ref{v1}) that
\begin{align}
\label{v22}
         &\var \left(\sum_{\floor{N/M_N}+1}^{\infty}\hat{f}(k)\*|t_{N,k}|^2\right)   =4\*\sum_{N/M_N<s\leq N-1}s^2|\hat{f}(s)|^2 + 4\*
(N^2-N)\sum_{N\leq s} |\hat{f}(s)|^2 \\
            &-4\sum_{\substack{N/M_N\leq s,t \\ 1\leq |s-t|\leq N-1\\ N\leq \max(s,t)}}(N-|s-t|)\hat{f}(s)\hat{f}(t)\hspace{2mm}
            -4\sum_{\substack{N/M_N\leq s,t\leq N-1\\N+1\leq s+t}}((s+t)-N) \hat{f}(s)\hat{f}(t). \nonumber
        \end{align}
The first term on the r.h.s. of (\ref{v22})  is equal to $2\*(V_N-V_{\floor{N/M_N}})$ and is $o(V_N)$
by (\ref{Kara}). The last two terms on the r.h.s. of (\ref{v22}) are $o(V_N)$ by Lemma 2.1.
Finally, the second term on the r.h.s. of (\ref{v22}) is bounded from above by 
\begin{align}
 4\*N^2\*\sum_{N\leq s} |\hat{f}(s)|^2\leq 2\* \sum_{j=1}^{\infty} \frac{1}{j^2}\*(V_{(j+1)\*N}-V_{j\*N})=o(V_N), 
\end{align}
where the last estimate follows from (\ref{karamata}).
This completes the proof of Lemma 2.2.
\end{pfo}

To finish the proof of the theorem, we have to show that
the first term on the r.h.s. of (\ref{aaa}) converges in distribution to a standard Gaussian random variable.  To do this, we first show that the first $\floor{M_N/2}$ moments of
\begin{align}
\label{part1}
\frac{2}{\sqrt{V_N}}\*\sum_{k=1}^{\floor{N/M_N}}\hat{f}(k)\*(|t_{N,k}|^2-\E |t_{N,k}|^2)
\end{align}
coincide with 
the first $\floor{M_N/2}$ moments of
\begin{align}
\label{expsum}
\frac{2}{\sqrt{V_N}}\*\sum_{k=1}^{\floor{N/M_N}}\hat{f}(k)k(\varphi_k-1),
\end{align}
where $\varphi_k, \ k\geq 1,$ are i.i.d. exponential random variables. 
\begin{lemma} 
Let $m$ be a positive integer such that $1\leq m<\frac{M_N}{2}.$ Then
\begin{align}
\label{momenty}
\E \left( \frac{2}{\sqrt{V_N}}\*\sum_{k=1}^{\floor{N/M_N}}\hat{f}(k)\*(|t_{N,k}|^2-\E |t_{N,k}|^2)\right)^m=\E \left(\frac{2}{\sqrt{V_N}}\*
\sum_{k=1}^{\floor{N/M_N}}\hat{f}(k)k(\varphi_k-1)\right)^m.
\end{align}
\end{lemma}
\begin{pfo}{\textit{Lemma 2.3}}
The formula (\ref{momenty}) follows from the identity
\begin{align}
\label{identity}
\E \prod_{i=1}^l |t_{N,k_i}|^2 =\E \prod_{i=1}^l k_i\*\varphi_{k_i},
\end{align}
provided 
\begin{align}
\label{condition}
2\*\sum_{i=1}^l |k_i|\leq N, \ \ 0<k_1, \ldots, k_l.
\end{align}
(\ref{condition}) was established in \cite{DS} and \cite{sasha} (see also \cite{jm}) where it was shown that a large number of the joint moments/cumulants of 
\[
k^{-1/2}\*t_{N,k}=k^{-1/2}\*Tr (U^k), \ \ k\geq 1,
\]
coincide with
the corresponding joint moments/cumulants of a sequence of i.i.d. standard complex Gaussian random variables. Namely, if we denote by $\kappa(t_{N,k_1}, \ldots, t_{N,k_n})$ the joint cumulant of
$\{t_{N,k_j}, \ 1\leq j\leq n\}$  then
\begin{align}
\label{kumulyanty}
\kappa(t_{N,k_1}, \ldots, t_{N,k_n})=0
\end{align}
if at one of the following two conditions is satisfied:
\begin{align*}
& (i) \ \ n\geq 1, \ \ \text{and} \ \sum_{j=1}^n k_j \neq 0 \\
& (ii) \ \ n>2, \ \ \sum_{j=1}^n k_j =0, \ \ \text{and} \ \sum_{j=1}^n |k_j|\leq N.
\end{align*}
In addition, $\kappa(t_{N,k}, t_{N,-k})=\min(|k|, N).$ We refer the reader to Lemma 5.2 of \cite{paper} for the details. 
Taking into account that the absolute value squared of a standard complex Gaussian random variable is distributed according to the exponential law we obtain (\ref{identity}-\ref{condition}).
This finishes the proof of Lemma 2.3.
\end{pfo}
Now, the proof of Theorem 1.2 immediately follows from the following standard lemma.
\begin{lemma}
For any $t\in \R$
\begin{align}
\E \exp\left( \frac{t}{\sqrt{\sum_{k=1}^N a_k^2}}\* \sum_{k=1}^{N} a_k\*(\varphi_k-1) \right) \to \exp(t^2/2)
\end{align}
as $N\to \infty,$ provided 
\begin{align}
\label{koef}
\sum_{k=1}^{\infty} a_k^2 =\infty \ \ \text{and} \ \ \max_{1\leq k\leq N} |a_k|=o\left(\sqrt{\sum_{k=1}^N a_k^2}\right).
\end{align}
\end{lemma}
\begin{pfo}{\textit{Lemma 2.4}}
The result of Lemma 2.4 follows from standard direct computations using independence of $\varphi_k$'s.
\end{pfo}

Setting $a_k=\hat{f}(k) \*k, \ \ k\geq 1$ and applying Lemma 2.4, we obtain the convergence of the exponential moment (and, hence, all moments) of $\sqrt{\frac{2}{V_N}}\*
\sum_{k=1}^{\floor{N/M_N}}\hat{f}(k)k(\varphi_k-1)$ to that of a standard real Gaussian random variable. Thus, both 
$$\sqrt{\frac{2}{V_N}}\*\sum_{k=1}^{\floor{N/M_N}}\hat{f}(k)k(\varphi_k-1)$$
and
$$\sqrt{\frac{2}{V_N}}\*\sum_{k=1}^{\floor{N/M_N}}\hat{f}(k)\*(|t_{N,k}|^2-\E |t_{N,k}|^2
$$
converge in distribution to a standard real Gaussian random variable. Since the second term in (\ref{aaa}) goes to $0$ in $L^2$ we conclude that
Theorem 1.2 is proven.
\end{pfo}

\end{document}